\documentclass[11pt]{amsart}
\usepackage{latexsym,amssymb,amsmath}
\usepackage{youngtab}
\input epsf
\def\boxit#1{\vbox{\hrule height1pt\hbox{\vrule width1pt\kern3pt
  \vbox{\kern3pt#1\kern3pt}\kern3pt\vrule width1pt}\hrule height1pt}}

\def\ta{\tilde{\a}}
\def\oa{\a_0}

\def\BC{\mathbb C}
\def\BA{\mathbb A}\def\BH{\mathbb H}
\def\BP{\mathbb P}
\def\pp#1{\mathbb P^{#1}}
\def\fz{\mathfrak z}

\def\fgl{\mathfrak g\mathfrak l}

\def\pp#1{{\mathbb P}^{#1}}
\def\tdim{{\rm dim}}

\def\cD{\mathcal D}

\def\trace{{\rm trace}}

\def\cJ{{\mathcal J}}

\def\CC{\mathbb C}

\def\HH{\mathbb H}

\def\OO{\mathbb O}

\def\ZZ{\mathbb Z}

\def\11{\mathbf 1}
\def\PP{\mathbb P}

\def\fh{{\mathfrak h}}

\def\fsl{{\mathfrak {sl}}}
\def\fsp{{\mathfrak {sp}}}

\def\fso{{\mathfrak {so}}}
\def\fe{{\mathfrak e}}
\def\fu{{\mathfrak u}}
\def\ff{{\mathfrak f}}

\def\fz{{\mathfrak z}}

\def\fg{{\mathfrak g}}
\def\fn{{\mathfrak n}}

\def\fk{{\mathfrak k}}

\def\a{\alpha}
\def\ta{\tilde{\alpha}}
\def\o{\omega}

\def\b{\beta}
\def\g{\gamma}
\def\s{\sigma}

\def\d{\delta}
\def\th{\theta}

\def\up#1{{}^{({#1})}}
\def\e{\varepsilon}
\def\ot{{\mathord{\,\otimes }\,}}
\def\op{{\mathord{\,\oplus }\,}}
\def\otc{{\mathord{\otimes\cdots\otimes}\;}}

\def\ra{{\mathord{\;\rightarrow\;}}}

\def\dim{{\rm dim}\;}

\newcommand\rem{{\medskip\noindent {\em Remark}.}\hspace{2mm}}
\newcommand\exam   {{\medskip\noindent {\em Example}.}\hspace{2mm}}

\newtheorem{theo}{Theorem}[section]

\newtheorem{lemm}[theo]{Lemma}
\newtheorem{prop}[theo]{Proposition}

\newtheorem{theorem}{Theorem}[section]
\newtheorem{proposition}[theorem]{Proposition}
\newtheorem{lemma}[theorem]{Lemma}
\newtheorem{corollary}[theorem]{Corollary}

\theoremstyle{definition}

\theoremstyle{remark}
\newtheorem{remark}[theorem]{Remark}

\begin{document}

\title{A Universal dimension formula \linebreak  for complex simple Lie algebras}
\author{J.M. Landsberg${}^1$ and L. Manivel}
\begin{abstract}
We present a universal formula for the dimension of the 
Cartan powers of the adjoint representation of a complex 
simple Lie algebra (i.e., a universal formula for the Hilbert 
functions of homogeneous complex contact manifolds), as well 
as several other universal formulas. These formulas generalize 
formulas of Vogel and Deligne and are given in terms of rational 
functions where both the numerator and denominator decompose into
 products of linear factors with integer coefficients. We   
discuss   consequences of the formulas including a relation with Scorza varieties.
\end{abstract}
 
\footnote{Supported by NSF grant DMS-0305829}
\maketitle

 \section{Statement of the main result}

Vogel \cite{vog1} defined a tensor category $\cD'$ intended to be a model for
a universal simple Lie algebra. His motivation came from knot theory,
as $\cD'$ was designed to surject onto the category of Vassiliev
invariants. While Vogel's work remains unfinished (and unpublished),
it already has consequences for representation theory.

Let $\fg$ be a complex simple Lie algebra. Vogel  derived 
a universal decomposition of $S^2\fg$ into   (possibly virtual)
Casimir eigenspaces, $S^2\fg = \BC\op Y_2(\a)\op Y_2(\b)\op Y_2(\g)$
which turns out to be a decomposition into
irreducible modules.  If we let $2t$ denote
the Casimir eigenvalue of the adjoint representation (with respect to some 
invariant quadratic form), these modules
respectively have Casimir eigenvalues $4t-2\a,4t-2\b,4t-2\g$, which
we may take as the definitions of $\a,\b,\g$.
Vogel showed   that $t=\a+\b+\g$.  He then went on
to find             Casimir eigenspaces $Y_3(\a),Y_3(\b),Y_3(\g)\subset S^3\fg$
with eigenvalues $6t-6\a,6t-6\b,6t-6\g$
(which again turn out to be irreducible), and computed 
their dimensions through difficult diagrammatic computations   and  
the help of Maple \cite{vog1}:

$${\rm dim}\, \fg  =  \frac{(\a-2t)(\b-2t)(\g-2t)}{\a\b\g}, $$
$${\rm dim}\, Y_2(\a)  =  -\frac{t(\b-2t)(\g-2t)(\b+t)(\g+t)(3\a-2t)}{\a^2\b\g(\a-\b)(\a-\g)}.$$ 
$${\rm dim}\, Y_3 (\a) =  -\frac{t(\a-2t)(\b-2t)(\g-2t)(\b+t)(\g+t)(t+\b-\a)(t+\g-\a)(5\a-2t)}
{\a^3\b\g(\a-\b)(\a-\g)(2\a-\b)(2\a-\g)},$$ 

\smallskip\noindent
and the formulas for $Y_2(\b),Y_2(\g)$ and $Y_3(\b),Y_3(\g)$ are obtained by permuting $\a$, $\b$, $\g$. 
These formulas suggest a completely different perspective from the usual description 
of the simple Lie algebras by their root systems and the Weyl dimension formula
that can be deduced for each particular simple Lie algebra. 
The work of Vogel raises many questions. In particular, what  remains 
of these formulas when we go to higher symmetric powers? If such
formulas do exist 
in general, do we need to go to higher and higher algebraic extensions
to state them, 
as Vogel suggests? Vogel describes  modules in the third tensor power
of the adjoint representation that require an algebraic extension
for their dimension formulas. 
 
For the exceptional series of simple Lie algebras, explicit 
computations of Deligne, Cohen and de Man showed that the decompositions
of the tensor 
powers are well-behaved up to degree 4, after whichÂ Â  modulesÂ Â  appear
whose dimensions
are {\it not} given by  rational functions whose numerator and
denominator are products of linear factors with integer coeffients 
(see \cite{del2,LMadv} for proofs of 
such types of formulas). In both the work of Vogel and Deligne et. al.,
problems arise when
there are different irreducible modules appearing in a Schur component
with the same Casimir eigenvalue. 

In this paper we   show that  some of the phenomena
observed by Vogel and Deligne do persist in all degrees.
Let  $\oa$ denote the highest root of $\fg$, once we have fixed a Cartan subalgebra 
and a set of positive roots.

\begin{theorem}
Use Vogel's parameters $\a,\b,\g$ as above. The $k$-th symmetric power
of $\fg$ contains three (virtual) modules
$Y_k(\a),Y_k(\b),Y_k(\g)$ with Casimir eigenvalues
$2kt-(k^2-k)\a,2kt-(k^2-k)\b,2kt-(k^2-k)\g$.
Using binomial coefficients defined by
$\binom{y+x}{y} = (1+x)\cdots (y+x)/y!$, we have:
$$\tdim\, Y_k(\a)=
\frac{t-(k-\frac 12)\a}{ t+\frac{\a}2}
\frac{\binom{-\frac{2t}{\a}-2+k }{ k}
\binom {\frac{  \b-2t }\a -1+k }{k }
\binom {\frac{  \g-2t }\a -1+k }{k } 
}
{
\binom{-\frac {\b}{\a} -1+k }{k }
\binom{-\frac {\g}{\a} -1+k }{k },
}$$
and   $\tdim\, Y_k(\b),\, \tdim\, Y_k(\g)$ are obtained by exchanging the
role  of $\a$  with $\b,\,\g$ respectively.
\end{theorem}

The modules $Y_k(\b),Y_k(\g)$ are described in \S \ref{whatare}.
For $Y_k(\a)$, we have the following refinement: 

\begin{theorem}
Parametrize the complex simple Lie algebras as follows:
$$\begin{array}{cccll}
Series & Lie\;algebra & \a & \b & \g \\
 & & & & \\
SP &   \fsp_{2n} & -2 & 1 &  n+2 \\
SL & \fsl_n & -2 & 2 & n\\
SO & \fso_n  & -2 & 4 & n-4 \\
EXC    &   &-2 & a+4 & 2a+4\\
  & \fsl_3 &-2 & 3 & 2\\
&\fg_2 & -2 & 10/3 & 8/3 \\
&\fso_8 & -2 & 4 & 4\\
 & \ff_4  & -2 & 5 & 6 \\
 & \fe_6  & -2 & 6 & 8 \\
 & \fe_7  & -2 & 8 & 12 \\
 & \fe_8  & -2 & 12  & 20 \\
F3_r  &  & -2 & a  & a(r-2)+4 
\end{array}$$
Then
   $Y_k(\a)$ is the $k$-th Cartan power $\fg^{(k)}$ of $\fg$ (the module
with highest weight $k\oa$) and
$$
\tdim\,\fg\up k = \frac{ \b+\g-3+2 k}{\b+\g-3}
\frac{ \binom{\b+  \g/2-3+k}{k} \binom{\g+  \b/2-3+k}{k}
 \binom{\b+\g-4+k}{k}}{\binom{-1+  \b/2+k}{k} \binom{-1+  \g/2+k}{k}}.
$$
\end{theorem}

In the exceptional series $EXC$, we have $a=-1, -2/3,0,1,2,4,8$.
Here $F3_r$ denotes the two-parameter series of Lie algebras in
the generalized third row of Freudenthal's magic chart,
$\fg_r(\BH,\BA)$ with $a= 1,2,4$ and $r\geq 3$,
which contains $ \fsp_{2r},\fsl_{2r},\fso_{4r}$, and $\fe_7$
when $r=3$ \cite{LMadv}. We call   $F3_3$ the {\it subexceptional
series}. 

The parameters $(\a,\b,\g)$ may be thought of as defining a
point in $\pp 2/\mathfrak{S}_3$ which we refer to as {\it Vogel's plane}.
We say a collection of 
points lie on a {\it line} in Vogel's plane if some lift of them
to $\pp 2$ is a colinear set of points. The classical series
$\fsl$, $\fso$, $\fsp$ all lie on lines by the description above
and one can even make the points of $\fso$ and $\fsp$ lie on the same
line. The algebras in the exceptional series all lie on a line, as do
the algebras in each of the generalized third rows of Freudenthal's magic 
chart.
Through
each classical simple Lie algebra there are an infinite number
of lines with at least three points.
Distinguished among these are the 
above-mentioned lines. For each
of these, there are natural inclusions of the 
Lie algebras as one travels north-east along the line.

\bigskip

$$ {\epsfxsize=5in\epsfbox{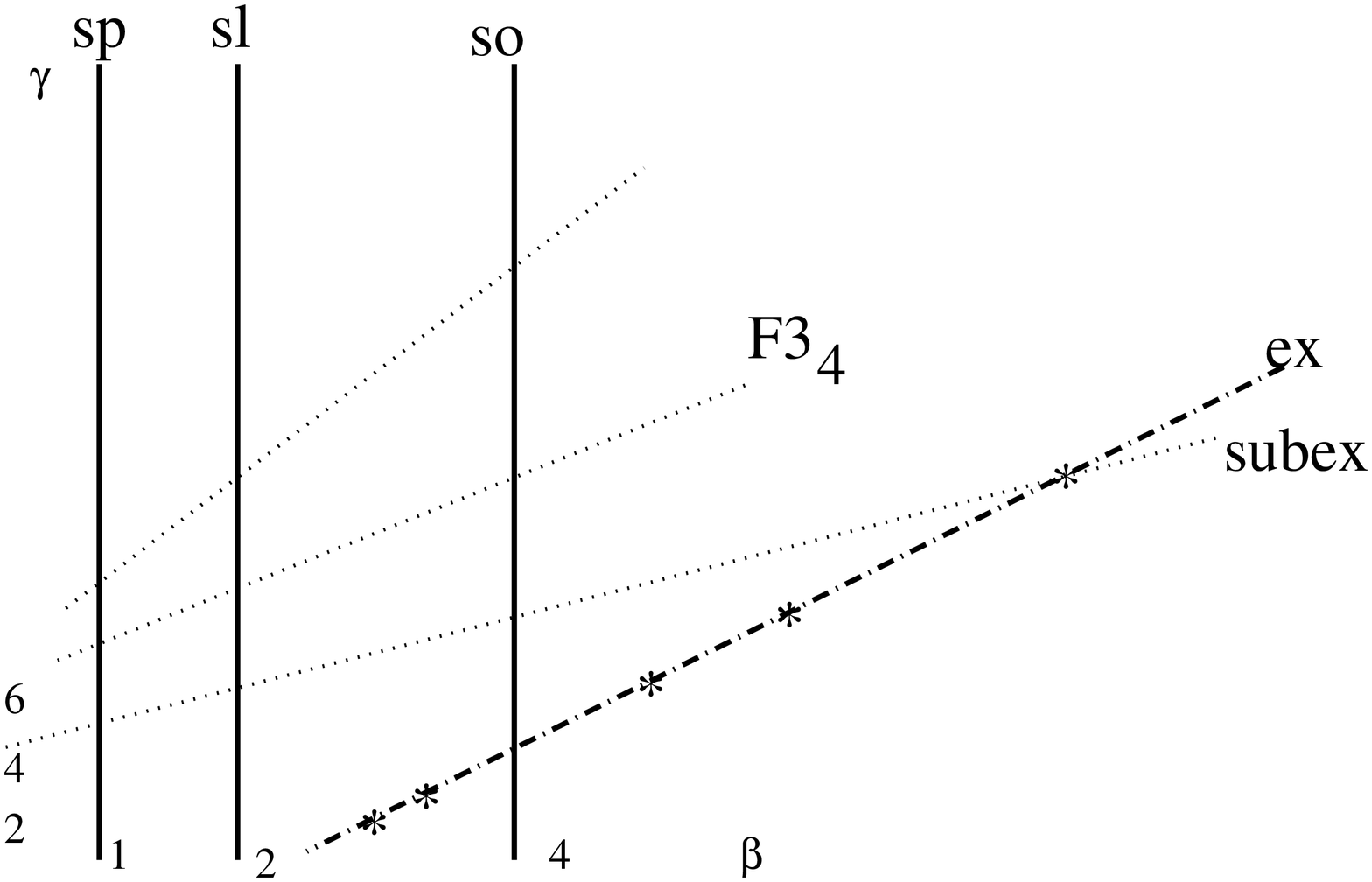}} $$

\centerline{{\sc Vogel's\; plane}}

Dotted diagonal lines correspond to $F3_r$.

\bigskip

\begin{remark}
The reason $\fso$ and $\fsp$ are split into two different lines is
we require $Y_k(\a)$ to be
the Cartan powers of the adjoint representation.
  The formula  for $\tdim \fg\up k$ applied to the $\fsp$ series
situated as $(-2,4,-2n)$ yields the dimensions of the modules $Y_k(\g)$.
\end{remark}

\medskip

\noindent{\bf Overview.}
In \S 2-4 we prove the main result, which is based on a careful 
analysis of the five step grading of a simple Lie algebra defined by a highest root. 
In \S 5 we show how this relates to other   
  $\ZZ$ and $\ZZ_2$-gradings and give a dimension formula
  for $\tdim \fg\up kY_2(\b)\up l$.  In  \S 6 we describe the modules $Y_k(\b),Y_k(\g)$ explicitly. 
We show that the highest weight of $Y_k(\b)$ is the sum of $k$ orthogonal long roots,
and give geometric interpretations
of them related to Scorza varieties. We conclude with an infinite
series of dimension formulas for the Cartan powers of the $Y_k(\b)$.
These formulas show that the modules $Y_k(\a)$ and $Y_k(\b)$ should be 
considered as universal in a very strong sense. Giving a 
precise meaning to that last sentence is an interesting  open problem. 

\medskip

\noindent {\bf Further questions and comments.}
Remarkably, the numbers $\b$ and $\g$ also appear  in \cite{kostant} 
in connection with the McKay correspondence. 
The numbers $h ,h' $ are exponents of $\fg$. For $\fg$ simply laced,
they coincide with the intermediate exponents of  the functions
$z(t)$ in \cite{kostant} having a linear factor. Why?

\smallskip
The formulas above, in addition to having zeros and poles, have indeterminacy
loci. For example,
the point corresponding to $\fso_8$
is in the indeterminacy locus of
  $\tdim Y_2(\b)$. 
For $\fso_8$,   $Y_2(\b)\op Y_2(\g)$ is the sum of the three isomorphic
$35$ dimensional representations  $2\o_1,2\o_3,2\o_4$. We obtain
$\tdim Y_2(\b)= 105$ (and $\tdim Y_2(\g)=0$)
 when considering $\fso_8$ as a member of the exceptional
series and $\tdim Y_2(\b)=70$ (and $\tdim Y_2(\g)=35$) when
considering it as an element of the orthogonal series. 
The same phenomenon occurs  for $\fsl_2$ which is also in the indeterminacy
loci of $\tdim Y_2(\b),\tdim Y_2(\g)$. While these remarks apply already to Vogel's
results (although we are unaware of them being pointed out before)   with the increasing
number of points in the indeterminacy loci as $k$ becomes large, it might be
interesting to address this issue in more detail.

\smallskip
We   remark that, for $k$ sufficiently negative the formulas above make sense
and give rise to dimensions of virtual modules. For example, in the
exceptional and subexceptional series, if one sets $K= 2t/\a +1-k$
then $\tdim Y_K(\a)= -\tdim Y_k(\a)$ and the dimensions for $k$ between $-1$
and $ 2t/\a$ are zero. Similar phenomena occur for the classical series.

\smallskip
Viewing the same equations with a different perspective,
we mention the work of Cvitanovic \cite{C76,cvit},
El Houari \cite{elh1,elh2} and Angelopolous \cite{angel}
which preceeded the work of Vogel and Deligne. Their works
contain calculations similar to Vogel's, but with a different goal: they
use that fact that dimensions of vector spaces are integers to classify
complex simple Lie algebras, and to
organize them into series, using Casimirs and
invariants of the symmetric algebra  to obtain diophantine equations.

\smallskip
If one restricts to the exceptional line, Cohen and deMan have observed
that (just using a finite number of dimension fomulas), the only
value of $a$ nontrivially yielding non-negative integers is, with our parametrization,
$a=6$. We account for this in \cite{LMsex} with a  Lie algebra 
which is intermediate between $\fe_7$ and $\fe_8$. It is 
an exceptional analogue of the {\it odd symplectic Lie algebras}. 
In fact the odd symplectic groups appear
to satisfy the formulas above 
when one allows $\g$ to be a half-integer
in the symplectic line. What other parameter values
yield integers in all the formulas?
Do the intermediate Lie algebras considered 
in \cite{LMsex} belong in Vogel's plane?  
\medskip

\smallskip
\noindent{\it Acknowledgements}: We thank Professor Deligne for his
comments on an earlier draft and Professor Gelfand for encouraging us
to continue this project after we had given up hope of finding
the main formula.

\section{The role of the principal $\fsl_2$.}

\subsection{How to use the Weyl   dimension formula}
A vector $X_{\oa}\in\fg_{\oa}$ belongs to
the minimal (nontrivial) nilpotent orbit in $\fg$. 
We can choose $X_{-\oa}\subset\fg_{-\oa}$ such that
$$(X_{\oa}, X_{-\oa}, H_{\oa}=[X_{\oa}, X_{-\oa}])
$$ 
is a $\fsl_2$-triple in $\fg$, 
generating a subalgebra of $\fg$ which we denote by $\fsl_2^*$. This is the {\it 
principal} $\fsl_2$. 
The semi-simple element $H_{\oa}$ defines a grading on $\fg$ according to the 
eigenvalues of $ad(H_{\oa})$:
$$\fg=\fg_{-2}\op  \fg_{-1}\op  \fg_0\op  \fg_1\op  \fg_2.$$

\smallskip\noindent The line $\fg_{2}$ (resp. $\fg_{-2}$) 
is generated by   $X_{\oa}$ (resp.  $X_{-\oa}$).
The subalgebra $\fg_0$ is reductive, and splits into the sum of the line generated 
by $H_{\oa}$  and the centralizer $\fh$ of the $\fsl_2$-triple. The $\fh$-module
$\fg_1$ is the sum of the root spaces $\fg_{\b}$, where $\b$ belongs to the set 
$\Phi_1$ of positive roots such that $\b(H_{\oa})=1$. Its dimension is twice the
dual Coxeter number of $\fg$, minus four \cite{knop}.

Let $\rho$ denote the half-sum of the positive roots. By the Weyl  dimension formula,

$$\dim\,\fg^{(k)}=\frac{(\rho + k\oa,\oa)}{(\rho,\oa)}
\prod_{\b\in\Phi_1}\frac{(\rho + k\oa,\b)}{(\rho,\b)}.$$

We thus need to analyze the distribution of the values of $(\rho,\b)$ for $\b\in \Phi_1$.

\subsection{The $\ZZ_2$-grading}
To do this, we slightly modify our grading of $\fg$. Let $V=\fg_1$ considered as   an irreducible 
$\fh\times \Gamma$-module, where $\Gamma$ is the automorphism
group of the Dynkin diagram of $\fg$. As an $\fh$-module, $V$ is irreducible
except in the case $\fg=\fsl_n$ where $\fh=\fgl_{n-2}$
and $V=V_{\o_1}\op V_{\o_{n-3}}$ as an $\fh$-module.
 
The space $V$ is endowed with a natural 
symplectic form $\o$ defined, up to scale, by the Lie bracket $\fg_1\times\fg_1\ra\fg_2$. 
Thus for each root $\b\in\Phi_1$, $\oa-\b$ is again a root in $\Phi_1$. 
Consider  $U=\fg_{\b-\oa}\op\fg_{\b}\subset\fg_{-1}\op\fg_1$. $U$ is stable under the adjoint 
action of $\fsl_2^*$, and is a copy of the natural two-dimensional $\fsl_2$-module. As a $\fsl_2^*\times\fh$-module,  
we thus get a $\ZZ_2$-grading of $\fg$ as 
$$\fg=\fg_{even}\op \fg_{odd}=\fsl_2^*\times\fh \op U\ot V.$$

 The Lie 
bracket defines an equivariant map 
$$\begin{array}{cccc}
\wedge^2(U\ot V) & = S^2U\ot \wedge^2V & \op & \wedge^2U\ot S^2V  \\
 & \qquad \quad\downarrow {\scriptstyle id\ot\o} & & \quad\downarrow {\scriptstyle \th}  \\
 & \quad\fsl_2^* & \op & \quad\fh. 
\end{array}$$
Here we use the natural identifications   $\fsl_2^*= S^2U$, and   $\wedge^2U=\CC$. 
Moreover, the fact that $\fh$ preserves the symplectic form $\o$ on $V$ implies that 
the image of $\fh$ in $End(V)\simeq V\ot V$ must be contained in $S^2V\simeq S^2V^*$. 
The map $\th$, up to scale, is dual to that inclusion. 

$$\begin{array}{ccccc}
Series & \fg & \fh & V & v \\
 & & & &  \\
SL & \fsl_n & \fgl_{n-2} & \CC^{n-2}\op (\CC^{n-2})^* & n-2\\
SO & \fso_n  & \fsl_2\times\fso_{n-4} & \CC^{2}\ot\CC^{n-4} & n-4 \\
SP &   \fsp_{2n} & \fsp_{2n-2} & \CC^{2n-2} & n-1 \\
EX   &\fg_2 & \fsl_2 & S^3\CC^2 & 2\\
&\fso_8 & \fsl_2\times\fsl_2\times\fsl_2  & \CC^{2}\ot\CC^{2}\ot\CC^{2} & 4 \\
 & \ff_4  & \fsp_6 & \wedge^{\langle 3\rangle}\CC^6 & 7 \\
 & \fe_6  & \fsl_6 & \wedge^{3}\CC^6 & 10 \\
 & \fe_7  & \fso_{12} & V_{\o_6}=\Delta_+ & 16 \\
 & \fe_8  & \fe_7  & V_{\o_7}=\fz_2(\OO) & 28 
\end{array}$$

\bigskip
Note that $\fsl_2^*\times\fh$ is a reductive subalgebra of maximal rank of $\fg$. 
We   choose   a Cartan subalgebra of $\fg$,  by taking the direct sum of $\CC H_{\oa}$, and 
  a Cartan subalgebra of $\fh$. The roots of $\fg$ will then be the root $\a_0=2\o_0$ 
of $\fsl_2^*$, the roots of $\fh$ and the weights of $U\ot V$, i.e., the sums 
$\pm\o_0+\mu$ with $\mu$ a weight of $V$. We can choose a set of positive roots 
of $\fh$, and if we choose the direction of $\a_0$ to be very positive,   
 the positive roots of $\fg$ will be $\a_0$, the positive roots of $\fh$, and the weights
$\o_0+\mu$ for $\mu$ any weight of $V$. Note that since $V$ is symplectic, the sum of these
weights must be zero.  Write $2v=\tdim V$, we have
$$2\rho=2\rho_{\fh}+(1+v)\a_0.$$
and
$$\Phi_1 =\{ \o_0+\mu \mid  \mu {\rm \ a \ weight \ of\ }V\}
$$

\smallskip
The set of simple roots of $\fg$ is easily described. If $\fg$ is not of type $A$, 
denote the highest weight of the irreducible $\fh$-module $V$ by $\chi$
so that its 
lowest weight is $-\chi$. For type
$A$ denote the highest weights by $\chi_1,\chi_2$, 
The simple roots of $\fg$ are
the simple roots of $\fh$ union $\o_0 -\chi$
($\o_0-\chi_1,\o_0-\chi_2$ for type $A$).
In particular the Dynkin diagram of $\fg$ is the diagram of $\fh$ with a vertex  
attached to the simple roots $\b$ of $\fg$ such that $(\chi,\beta)
\neq 0$, with the obvious analog attaching two verticies for type $A$.

\rem If we had chosen the directions of $\fh$ to be much more positive than that of 
$\fsl_2^*$, we would have obtained a different set of positive roots and, 
{\it except for} $\fg=\fg_2$,  
the highest root $\ta$ of $\fh$  would have been the highest root of $\fg$
(here we suppose that $\fh$ itself is simple; otherwise we can take the highest root
of any simple factor of $\fh$). We suppose in the sequel that 
we are not in type $\fg_2$: then
$\a_0$ and $\ta$, considered as roots of $\fg$, are both long.

\medskip

For type $C$, the root $\o_0+\chi$ is short and it is long in all other cases.
 This is because, except in type $A$ which can be checked seperately,
  $\a_0= k\o_i$ for  some fundamental weight $\o_i$, and
$\o_0+\chi$, being the second highest root, must equal $k\o_i-\a_i$. But then, 
$$(k\o_i-\a_i,k\o_i-\a_i)=(k\o_i,k\o_i)+(1-k)(\a_i,\a_i).$$
We conclude that $\o_0+\chi$ is long iff the adjoint representation is fundamental, 
i.e.,
 iff we are not in type $C$. Then
$$(\o_0+\chi,\o_0+\chi)=(\a_0,\a_0)=(\ta,\ta)=\frac{4}{3}(\chi,\chi).$$
Another interesting relation can be deduced from the fact that for any simple 
root $\a$ of $\fg$, we have $(2\rho,\a)=(\a,\a)$, since $\rho$ is the sum of the 
fundamental weights. Applying this to $\a=\o_0-\chi$, we get
$$(\chi+2\rho_h,\chi)=\frac{2v+1}{4}(\a_0,\a_0).$$

\medskip
Note that the scalar form here is the Killing form of $\fg$, more precisely the
dual of its restriction to the Cartan subalgebra. Restricted to the duals of the 
Cartan  subalgebras of $\fsl_2^*$ or $\fh$, we can compare it to their Killing forms.
Suppose that $\fh=\fh_1\times\cdots\times\fh_m$, and $V=V_1\otc V_m$ for some 
$\fh_i$-modules $V_i$. 

To simplify notation in the calculations that follow we use the normalization
that the Casimir eigenvalue of every simple Lie algebra is $1$, i.e., we use
for  invariant quadratic form
the Killing form $K(X,Y)=\trace (ad(X)\circ ad(Y))$.

Then  for $X\in\fsl_2^*$ and $Y\in\fh_i$, we have 
\begin{eqnarray}
\nonumber & \trace _{\fg}ad(X)^2   =   \trace _{\fsl_2^*}ad(X)^2+2v\trace_U X^2  =  
(1+\frac{v}{2})\trace _{\fsl_2^*}ad(X)^2, \\
\nonumber &
\trace _{\fg}ad(Y)^2  =    \trace_{\fh_i}ad(Y)^2+2\frac{\dim V}{\dim V_i}\trace_{V_i} Y^2  =  
(1+4ve_{V_i})\trace _{\fh}ad(Y)^2,
\end{eqnarray}
where $e_{V_i}$ is related to the Casimir eigenvalue $c_{V_i}$ of $V_i$ by the identity
$e_{V_i}=\frac{c_{V_i}}{\dim\fh_i}.$
Taking duals, we deduce that 
$$\begin{array}{rcl}
(\a_0,\a_0) & =&\frac{2}{v+2}(\a_0,\a_0)_{\fsl_2^*}=\frac{1}{v+2}, \\
 & & \\
(\ta,\ta)_{\fh}  &= & \frac{4}{3}(\chi,\chi)_{\fh}=\frac{1+4ve_V}{v+2}.
\end{array}$$

\smallskip
Note that the dual Coxeter number of a simple Lie algebra is 
the Casimir eigenvalue of the Lie algebra divided by the length
of the longest root. We conclude that the dual Coxeter 
number $ \check h$ of $\fg$ is $ \check h = v+2$, while the dual Coxeter number of $\fh$, 
which 
we denote by $h$, is equal to $\frac{1+4ve_V}{v+2}$. Remember that the normalization 
of the Killing form is such that $(\ta+2\rho_{\fh},\ta)_{\fh}=1$, so that 
$(2\rho_{\fh},\ta)_{\fh}=(h-1)(\ta,\ta)_{\fh}$, thus $(2\rho_{\fh},\ta)=(h-1)(\ta,\ta)$
as well. 
\medskip

\section{The Casimir eigenvalues of $S^2\fg$}

\subsection{A nontrivial component in the symmetric square of $\fg$}
Vogel proved that $S^2\fg$ can contain at most four Casimir eigenspaces 
  (allowing the possibility of zero, or even virtual eigenspaces).
Two irreducible components are obvious: the Cartan square, whose highest,
weight is $2\a_0$, and the trivial line generated by the Killing form. 
We identify, for $\fg$ not of type $A_1$ (i.e., $\fh\ne 0$), another component. 

\begin{proposition} \label{anotherone}
The symmetric square $S^2\fg$ has a component $Y_2(\b)$ of highest 
weight $\a_0+\ta$. 
\end{proposition}

\proof From our $\ZZ_2$-grading $\fg=\fsl_2\times\fh\op U\ot V$, we deduce that 
$$S^2\fg=S^2\fsl_2\op S^2\fh\op (S^2U\ot S^2V)\op (\wedge^2U\ot\wedge^2V)\op
(\fsl_2\ot\fh)\op (\fsl_2\ot U\ot V)\op (U\ot\fh\ot V).$$
All the weights here are of the form $k\o_0+\mu$ for $\mu$ in the weight lattice of $\fh$, 
and we will call the integer $k$ the {\it level} of the weight. The maximal level
is four, and the unique weight of level four is $2\a_0$, the highest weight of 
$S^2\fsl_2$. The corresponding weight space, of dimension one, generates the Cartan
square of $\fg$.  

We will check that once we have suppressed the weights of the Cartan square  with
their multiplicities, the highest remaining weight is $\a_0+\ta$, which has level 
two. 

At level three, we only get weights coming from $\fsl_2\ot U\ot V$. More precisely, 
let $e,f$ be a basis of $U$ diagonalizing our Cartan subalgebra, in such a way that 
the semi-simple element $H$ of our $\fsl_2$-triple has eigenvalues $1$ on $e$, $-1$ on 
$f$, while $X=f^*\ot e$ and $Y=e^*\ot f$. Then a weight vector of level three in 
$\fsl_2\ot U\ot V$ is of the form $X\ot e\ot v$, for some weight vector $v\in V$, and such
a weight vector is contained in $\fg.X^2\subset S^2\fg$, hence in the Cartan square
of $\fg$. It is equal, up to a nonzero constant, to $(f\ot v).X^2$. We conclude
that all weight vectors of level three belong to the Cartan square of $\fg$. 

We turn to level two. First observe that $\a_0+\ta$ has multiplicity two inside
$S^2\fg$. It is the highest weight of $\fsl_2\ot\fh$, which appears twice in the 
decomposition of $S^2\fg$ above: once  as such, and once in a slightly more hidden way, 
as a component of $S^2U\ot S^2V$. Indeed, recall that $S^2U$ and $\fsl_2$ are equal, 
and that we defined a nontrivial map $\th :S^2V\ra\fh$. We   check that $\a_0+\ta$ 
has only multiplicity one inside the Cartan square of $\fg$, and our claim will follow. 

Note that this Cartan square is $U(\fn_-)X^2$, where $\fn_-\subset\fg$ is the subalgebra
generated by the negative root spaces and $U(\fn_-)$ its universal envelopping algebra. 
As a vector space, this algebra is generated by monomials on vectors of negative weight, 
hence of nonpositive level. How can we go from $X^2$, which is of level four, to 
some vector of level two? We have to apply a vector of level $-2$, or twice a
vector of level $-1$. For the first case, the only possible vector is $Y$, which maps 
$X^2$ to $XH\in S^2\fsl_2$. For the second case, we first apply some vector $f\ot v$,
with $v\in V$: this takes $X^2$ to $X\ot (e\ot v)$, up to some constant. Then we 
apply another vector $f\ot v'$, and obtain, again up to some fixed constants, 
$$(e\ot v')(e\ot v)+X\ot\th (vv')+\o(v,v')XH.$$
The first component belongs to $S^2(U\ot V)$, the second one to $\fsl_2\ot\fh$, 
the third one to $S^2\fsl_2$. The contribution of the first component to $\fsl_2\ot\fh
\subset S^2(U\ot V)$ is $e^2\ot\th (vv')=X\ot\th (vv')$. We conclude that the Cartan
square of $\fg$ does not contain $\fsl_2\ot\fh\op\fsl_2\ot\fh \subset S^2U\ot S^2V
\op\fsl_2\ot\fh$, but meets it along some diagonal copy of $\fsl_2\ot\fh$. This 
implies our claim. \qed

\subsection{Interpretation of Vogel's parameters} 
It is now easy to compute the Casimir eigenvalues of our two 
nontrivial components of $S^2\fg$:
$$\begin{array}{rcl}
C_{Y_2(\a)}& = & (2\a_0+2\rho,2\a_0)=2\frac{v+3}{v+2}, \\
C_{Y'_2}& = & (\a_0+\ta+2\rho,\a_0+\ta)=\frac{v+h+2}{v+2}.
\end{array}$$

\begin{corollary} Let $h'=v-h $. Normalize Vogel's parameters for 
$\fg\ne\fg_2$, such that $\a=-2$. Then 
$$\quad \b=h'+2, \quad \g=h+2, \quad t=v+2=\check h.$$
\end{corollary}

\proof Vogel's parameters are defined by the fact that, 
with respect to an   invariant quadratic form on $\fg$, 
the Casimir eigenvalue of $\fg$ is $2t$,   the nonzero Casimir eigenvalues
of $S^2\fg$   are $2(2t-\a), 2(2t-\b), 2(2t-\g)$, and $t=\a+\b+\g$. 
We have been working with the Killing form, for which the Casimir eigenvalue of $\fg$ is 
$1$. 
Rescaling $t$ to be $v+2$ and plugging into the formulas for $C_{Y_2(\a)},C_{Y_2(\b)}$
we obtain the result.
\qed

\medskip\rem 
Note that this does not depend on the fact that $\fh$ is simple. 
If it is not simple, we can choose a highest root for any simple factor and
get  a corresponding component of $S^2\fg$, whose Casimir eigenvalue is given as before 
in terms of the dual Coxeter number of the chosen factor. This implies that we cannot
have more than two simple factors, and that when we have two, with dual Coxeter  
numbers $h_1$ and $h_2$, then $v=h_1+h_2$. This actually happens in type $B$ or $D$. 
(Beware that this should be understood up to the symmetry of the Dynkin diagram: 
in type $D_4$ we get three different components in $S^2\fg$, but they are permuted 
by the triality automorphisms and their sum must be considered as simple.)  

\bigskip
If $\fh$ is simple, the formula $(\chi+2\rho_{\fh},\chi)=\frac{2v+1}{4}(\ta,\ta)$
gives $c_V=\frac{2v+1}{4h}$, and we get 
$$\dim\fh =\frac{c_V}{e_V}=\frac{v(2v+1)}{h'+2}.$$

In general, Vogel's dimension formula is 
$$\dim\fg=d(h,h')=\frac{(h+h'+3)(2h+h'+2)(h+2h'+2)}{(h+2)(h'+2)}.$$

We have the following curious consequence. Parametrize
$\fg$ by $h$ and $h'$. We ask: What values of $h$ and $h'$
can give rise to a $\fg$ such that $\fh$ is simple and $V$ is
irreducible? In this case $\fh$ is parametrized by 
$h'$ and $h-h'$. Thus $d(h,h')=d(h',h-h')+3+4(h+h')$, which is equivalent 
to the identity $(h+1)(h-2h'+2)=0$. Thus such $\fg$ must be
in  the symplectic series $h=-1$, 
or the exceptional series $h=2h'-2$! 

\subsection{Interpretation of $h'$} 
Suppose that we are not in type $A$, so that the adjoint representation is 
supported on a fundamental weight $\o$. 
Let  $\a_{ad}=\o_0-\chi$ denote the corresponding simple root
dual to $\o$. Since the highest root 
$\ta$ of $\fh$ is not the highest root of $\fg$, one can find  a simple root $\a$ such that 
$\ta+\a$ is again a root, and the only possibility is $\a=\a_{ad}$. Thus $\ta+\a_{ad}$,
and by symmetry $\psi=\a_0-\ta-\a_{ad}$, both belong to $\Phi_1$. Suppose that $V=V_{\chi}$ is 
fundamental, and let $\a_{\chi}$ be the corresponding simple root.

\begin{prop} 
$\phi=\psi-\a_{ad}-\a_{\chi}$ is the highest root 
of $\fg$ orthogonal to $\a_0$ and $\ta$.
\end{prop} 

\proof We first prove that $\phi$ is a root. First note that 
$$(\psi, \a_{\chi})=(\chi-\ta,\a_{\chi})=(\a_{\chi},\a_{\chi})/2-(\ta,\a_{\chi}).$$
If we are not in type C, then $(\ta,\a_{\chi})=0$. Indeed, $\ta$ is a fundamental
weight and $\a_{\chi}$ is a simple root, so if this were non-zero we would get 
$\ta=\chi$, which cannot be since we know that $\chi$ is minuscule. 
Thus $(\psi, \a_{\chi})>0$ and $\psi'=\psi-\a_{\chi}$ is a root. Moreover,
$$(\psi, \a_{ad})=(\a_0,\a_0)/4-(\chi,\chi)+(\chi,\ta)=-(\a_0,\a_0)/2+(\ta,\ta)
\chi (H_{\ta})/2,$$
where $\chi (H_{\ta})$ is a positive integer (in fact equal to one, since we know that 
$V$ is minuscule). Thus $(\psi, \a_{ad})\ge 0$ and 
$(\psi', \a_{ad})\ge -(\a_{\chi},\a_{ad})=(\a_{\chi},\chi)>0$. We conclude that 
$\psi'-\a_{ad}=\phi$ is a root. 
  
Since $\phi=2\chi-\ta-\a_{\chi}$, it is clearly orthogonal to $\a_0$. Moreover, 
using again that $(\ta,\a_{\chi})=0$ if we are not in type C, we have 
$(\phi,\ta)=(2\chi-\ta,\ta)=0$, since we have just computed that $(\chi,\ta)=
(\ta,\ta)/2$. To conclude that $\phi$ is the highest root orthogonal to both 
$\oa$ and $\ta$, we use the following characterization of the highest root 
of a root system:

\begin{lemm}
The highest root of an irreducible root system is the only long root $\g$
such that $(\g,\a)\ge 0$ for any simple root $\a$. 
\end{lemm}

We apply this lemma to $\phi=2\chi-\ta-\a_{\chi}
=-\ta-\sum_{\a\ne\a_{\chi}}c_{\a,\a_{\chi}}\o_{\a}$, where $\a$ belongs
to the set of simple roots and $c_{\a,\a_{\chi}}$ is the corresponding Cartan
integer. Since $\a\ne\a_{\chi}$, we know that  $c_{\a,\a_{\chi}}\le 0$, 
hence $(\phi,\a)\ge 0$ for every simple root $\a\ne\a_{ad}$. 

It  remains to check that $\phi$ is long. Remember that $(\chi,\chi)=\frac{3}{4}
(\ta,\ta)$, that $(\ta,\a_{\chi})=0$ and $(\ta,\chi)=(\ta,\ta)/2$. We compute that 
$(\phi,\phi)=2(\ta,\ta)-(\a_{\chi},\a_{\chi})\ge (\ta,\ta)$. Therefore $\phi$
is long (and we must have equality, so that $\a_{\chi}$ is also long). \qed

\medskip Note that, say in the simply laced case,  $$(2\rho,\phi)=(2\rho,\a_0)-(2\rho,\ta)
-3(\a_0,\a_0)=(v+1-h+1-3)(\a_0,\a_0)=(h'-1)(\a_0,\a_0).$$ We have therefore isolated 
three roots $\a_0, \ta, \phi$ of heights $v+1, h-1, h'-1$ respectively.

\section{Proof of the main result}

\subsection{The weights of $V$ and their heights} 
Our next observation concerns the distribution of the rational numbers $(\rho_{\fh},\mu)$, 
when $\mu$ describes the set of weights of $V$. A natural scale for these numbers is the 
length $(\ta,\ta)$ of the long roots. We denote by $S_p$ the string of numbers 
$(p/2-x)(\ta,\ta)/2$, for $x=0,1,\ldots ,p$. 

\begin{proposition} \label{atmosttwolem}
The values $(\rho_{\fh},\mu)$, for $\mu$ a weight of $V$, can be arranged into the union of the 
three strings $S_{ {v-1} }$, $S_{ {h-1} }$, $S_{ {h'-1} }$.
\end{proposition}

Our main theorem easily follows from this fact: a set of weights $\mu$ in $V$ contributing 
to a string $S_p$ of values of $(\rho_{\fh},\mu)$, gives a set of roots $\b=\o_0+\mu$ in 
$\Phi_1$, with 
$$\begin{array}{rcl}
(\ta,\b) & = & (\a_0,\a_0)/2, \\
(\rho,\b) & = & \frac{1+v}{2}(\a_0,\a_0)/2+(\rho_{\fh},\mu)=(\frac{1+v}{2}+
p/2-x)(\a_0,\a_0)/2, \quad 0\le x\le p. 
\end{array}$$
The contribution of this subset of $\Phi_1$ to the Weyl dimension formula is therefore 
$$C_p=\prod_{\b}\frac{(\rho + k\ta,\b)}{(\rho,\b)}=\prod_{x=0}^p
\frac{\frac{1+v+p}{2}-x+k}{\frac{1+v+p}{2}-x}=\frac{\binom{\frac{v+p+1}{2}+k}{k}}
{\binom{\frac{v-p-1}{2}+k}{k}}.$$

\medskip Our proof of the Proposition is a case by case check. One has to be careful 
about the case where $h'\le 0$, since the string $S_{\frac{h'-1}2}$ is no longer defined. 
Since for $p>0$, we have $C_{-p}=C_{p+1}^{-1}$, we should interpret $S_{-p}$ as 
{\it suppressing a string} $S_{p+1}$. We then easily check the rather surprising 
fact that, interpreted that way, the Proposition also holds for $h'<0$.

\subsection{Relation with Knop's construction of simple singularities}
F. Holweck observed that the fact that we can arrange the values of $(\rho,\b)$, for
$\b\in\Phi_1$, in no more than three strings, has a curious relation with the work
of F. Knop on simple singularities. 
Knop proved \cite{knop} that if $Y^{\perp}\subset \BP\fg$ is
a hyperplane Killing orthogonal to a regular nilpotent element $Y\in\fg$, 
the intersection of this hyperplane with the adjoint variety $X_{ad}\subset
\PP\fg$ (the projectivization of the minimal nontrivial nilpotent orbit),
has an isolated singularity which is simple, of type given by the subdiagram of 
the Dynkin diagram of $\fg$  obtained from   the long simple roots. 

We can choose $Y=\sum_{\a\in\Delta}X_{\a}$, where $\Delta$ denotes the set of 
simple roots and $X_{\a}$ is a generator of the root space $\fg_{\a}$. 
The orthogonal hyperplane contains the lowest root space $\fg_{-\ta}\in X_{ad}$.
Let $P$ denote the parabolic subgroup of the adjoint group of $\fg$, which 
stabilizes $\fg_{\ta}$, and let $U$ denote its unipotent radical. Being unipotent,
$U$ can be identified, through the exponential map, 
with its algebra $\fu$, a basis of which is given by the 
roots spaces $\fg_{\b}$ with $\b\in\Phi_1\cup\{\a_0\}$. The scalar product 
with $Y$ defines on the Lie algebra $\fu$ the function $f(X)=K(Y,exp(X)X_{\ta})$.
The quadratic part of this function is
$$q(X,X')  = \frac{1}{2} K(Y,ad(X')ad(X)X_{\ta}) =
\frac{1}{2}K(ad(X)Y,ad(X')X_{\ta}).$$
The kernel of this quadratic form thus contains the 
kernel of the map $X\mapsto [Y,X]$, $X \in \fu$. 

Suppose for simplicity that $\fg$ is simply laced. Let 
$$\fg_l=\bigoplus_{\g\in\Phi_1\cup\{\a_0\},\;
(\rho,\g)=l}\fg_{\g}.$$ 
Then $ad(Y)$ maps $\fg_l$ to $\fg_{l+(\ta,\ta)/2}$. In particular,
the kernel of $ad(Y)_{|\fg_l}$ has dimension at least $\dim \fg_l-\dim \fg_{l+(\ta,\ta)/2}$.
Since $\fg_v=\fg_{\a_0}$ is one-dimensional, we deduce that for all $l$, 
$$\dim\fg_l\le 1+{\rm corank}(q).$$
The maximal dimension of $\fg_l$ is the minimal number of strings we need to 
arrange the values of $(\rho,\b)$ for $\b\in\Phi_1$. This number is bounded by
three because, since $f$ defines a simple singularity, the corank of its quadratic
part must be  at most two. 

Note that in Knop's work there is no direct proof of this fact. It follows
from a numerical criterion and a trick attributed to K. Saito (\cite{knop}, Lemma 1.5). 
\medskip

\section{Gradings}

The highest roots $\a_0$ and $\ta$ both induce $5$-gradings on $\fg$. 
Being orthogonal, they induce a double grading 
$$\fg_{ij}=\{X\in\fg, [H_{\a_0},X]=iX, [H_{\ta},X]=jX\}.$$

\begin{proposition} \label{prop51} With the normalization $t=\check h$,
for $\fg$ of rank at least three,
the dimensions of the components of this double grading are given by the 
following diamond:
$$\begin{array}{ccccc}
 & & 1 & & \\
 & \b & 2\g-8 & \b & \\
1 & 2\g-8 & \ast & 2\g-8 & 1 \\ 
 & \b & 2\g-8 & \b & \\
 & & 1 & & 
\end{array}$$
\end{proposition}
 
\proof Let $g_{ij}$ denote the dimension of $\fg_{ij}$. 
Since the dual Coxeter number of $\fg$ is $t $, 
the dimension of the positive part of the $5$-grading of $\fg$
is $2g_{11}+g_{01}+1=2t-3=2\b+2\g-7$. Since the dual Coxeter number of $\fh$ is $h $,  
the dimension of the  
positive part of the $5$-grading of $\fh$ is $g_{01}+1=2h-3=2\g-7$.
Hence the claim. \qed 

\begin{corollary} 
With the normalization $t=\check h$, the integer $\b$ is the number of roots $\theta$ in $\Phi_1$ 
such that $\ta+\theta$ is still a root. 
\end{corollary}

\proof Let $\theta\in \Phi_1$ be such that $\fg_{\theta}\subset\fg_{11}$. 
This means that $\theta(H_{\ta})=1$. Then $\theta^*=\a_0-\theta$ is also
a root, and $\theta^*(H_{\ta})=-1$, thus $s_{\ta}(\theta^*)=\ta+\theta^*$
is a root. Conversely, if $\ta+\theta^*$ is a root, $(\ta+\theta^*)(H_{\ta})
=2+\theta^*(H_{\ta})$, as well as $\theta^*(H_{\ta})$, belongs to $\{-2,-1,0,1,2\}$, 
hence $\theta^*(H_{\ta})=-1$ and we can recover $\theta\in\Phi_1$.
The $\fg_2$ case may be verified directly. \qed

\medskip Let  $\fg_{00}^*\subset\fg_{00} $
denote the common centralizer of
$X_{\a_0}$  and $X_{\ta}$. We have $\fg_{00}=\fg_{00}^*\op\CC H_{\a_0}\op\CC H_{\ta}$.

\begin{proposition}
$\fg_{11}$ is endowed with a $\fg_{00}^*$-invariant 
non-degenerate quadratic form.
\end{proposition}

\proof For $Y,Z\in\fg_{11}$, let 
$$Q(Y,Z)=K([X_{\ta},Y],[X_{-\a_0},Z]).$$
This bilinear form is obviously $\fg_{00}^*$-invariant. We check it is symmetric:
$$\begin{array}{rcl}
Q(Y,Z) & = & K(X_{\ta},[Y,[X_{-\a_0},Z]]) \\
 & = & K(X_{\ta},[Z,[X_{-\a_0},Y]])+ K(X_{\ta},[X_{-\a_0},[Y,Z]])\\
 & = & Q(Z,Y)+\Omega(Y,Z)K(X_{\ta},H_{\a_0}) \\
 & = & Q(Z,Y).
\end{array}$$
Recall that $K(\fg_{\b},\fg_{\g})\ne 0$ if and only if $\b+\g=0$. To prove
that $Q$ is non-degenerate, we must therefore check that for each root space 
$\fg_{\theta}$ in $\fg_{11}$, $\ta+\theta$ is a root -- this follows from 
the corollary above -- and $\a_0-(\ta+\theta)$ is also a root -- this follows 
from the fact that $\ta+\theta$ is in $\Phi_1$. \qed

\medskip We thus get an invariant map $\fg_{00}^*\rightarrow\fso_{\beta}$, 
which turns out to be {\it surjective}. We can thus 
write $\fg_{00}=\CC^2\times\fso_{\b}\times\fk$ for some reductive subalgebra
$\fk$ of $\fg$. Our double grading of $\fg$ takes the form:

$$\begin{array}{ccccc}
 & & \CC & & \\
 & \CC^{\b} & U^{2\g-8} & \CC^{\b} & \\
\CC & U^{2\g-8} &\CC^2\times\fso_{\b}\times\fk  & U^{2\g-8} & \CC \\ 
 & \CC^{\b} & U^{2\g-8} & \CC^{\b} & \\
 & & \CC & & 
\end{array}$$
Note that $U$ is a symplectic $\fk$-module. 
Now, consider the $5$-step simple grading that we obtain by taking diagonals. 
Since $\fso_{\b}\op \CC^{\b}\op\CC^{\b}\op\CC=\fso_{\b+2}$, we get 
$$\fg=\CC^{\b+2}\op V^{4\g-16}\op (\CC\times\fso_{\b+2}\times\fk)\op V^{4\g-16}
\op\CC^{\b+2}.$$
Remarkably, this induces a very simple $\ZZ_2$-grading 
$$\fg=(\fso_{\b+4}\times\fk)\op W^{8\g-32}.$$
The subalgebra $\fk$ and the module $W^{8\g-32}$ are given by the following table:

$$\begin{array}{ccccc}
\b & \fg & \fk & \fso_{\b+4} & W \\ 
1  & \fsp_{2n} & \fsp_{2n-4} & \fso_5=\fsp_4 & A_4\ot B_{2n-4}\\
2  & \fsl_{n} & \fgl_{n-4} & \fso_6=\fsl_4 & A_4^*\ot B_{n-4}\op A_4\ot B_{n-4}^*
\\
4  & \fso_{n} & \fso_{n-8} & \fso_8 & A_8\ot B_{n-8}\\
5 & \ff_4 &Â  &Â  \fso_9 & \Delta \\
6 & \fe_6 & \CC &Â  \fso_{10} & \Delta_+\op\Delta_- \\ 
8 & \fe_7 & \fsl_2 &Â  \fso_{12} & A_2\ot \Delta_+ \\ 
12 & \fe_8 &Â  &Â  \fso_{16} & \Delta_+ \\ 
a & F3G(a,r) &F3G(a,r-2)Â  &Â  \fso_{a+4} & \BC^{8a(r-2)} 
\end{array}$$

\subsection{More dimension formulas}
Recall from Proposition 3.1 that $Y_2(\b)$ has highest weight $\ta+\a_0$.
Therefore, the diamond of Proposition 5.1 also gives the number of
roots having a given scalar product with the highest weights of 
$\fg$ and  $Y'_2$.  From Corollary 4.1, we know the values 
of the $(\rho,\a)$ when $\a$ describes the roots in $\Phi_1$. 
Among these, the $\b$ positive roots from $\fg_{1,-1}$ (and $\fg_{1,1}$,
symetrically) contribute a string of length $\b-1$, the middle point 
having multiplicity two (just like the
weights of the natural representation of $\fso_{\b}$, in accordance
with Proposition 5.3). Finally, since $\fg_{0,1}$ is part of the 
5-step adjoint grading of $\fh$, its contribution can be 
described by Proposition 4.1 with the pair $h,h'$ changed into
$h',h-h'$. Using the Weyl  dimension formula, we get: 

\begin{theorem} 
The dimensions of the Cartan products of powers of $\fg$ and $Y'_2$ are 
given by the universal formula 
$$\tdim\,\fg\up kY'_2\up l  =  F(\b,\g,k,l)
A(\b,\g,k+l)B(\b,\g,l)C(\b,\g,k+2l)C(\b,\g,k-\g+3), $$
with
$$F(\b,\g,k,l)  =   \frac{(\b+\g-3+2k+2l)(\g-3+2l)(\b/2+\g-3+k+2l)(\b/2+k)}
{(\b+\g-3)(\g-3)(\b/2+\g-3)\b/2}, $$

\begin{eqnarray}
\nonumber
A(\b,\g,k) &=
&\frac{\binom{\b+\g/2-3+k}{k}\binom{\g+\b/2-4+k}{k}\binom{\g-3+k}{k}}
{\binom{\b/2+k}{k}\binom{-1+\b+k}{k} \binom{-1+\g/2+k}{k}}, \\
\nonumber
B(\b,\g,k) &=
&\frac{\binom{\g-\b/2-3+k}{k}\binom{\g/2+\b/2-3+k}{k}\binom{\g-4+k}{k}}
{\binom{\b/2-1+k}{k}\binom{\g/2-\b/2-1+k}{k}},\\
\nonumber
C(\b,\g,k) &= &\frac{\binom{\b+\g-4+k}{k}}{\binom{\g-3+k}{k}}.
\end{eqnarray}
\end{theorem}

\section{The modules $Y_k(\b)$ and $Y_k(\g)$}\label{whatare} 
 
For each simple Lie algebra $\fg$ we have obtained a general formula for
the dimension of its 
$k$-th Cartan power as a rational function of $\a,\b,\g$, symmetric with respect to
$\b$ and $\g$.  Following Vogel, the three numbers should play a completely
symmetric role, and by permutation we should get the dimensions of (virtual)
$\fg$-modules $Y_k(\b)$ and $Y_k(\g)$. We first check that this is indeed the case.
The formula predicts that these modules must be zero when $k$ becomes large, 
but an interesting pattern shows up in the classical cases.

\subsection{Identification}
The formula for the dimension of $Y'_k$ is 

$$\dim\,Y'_k = \frac{2\g-(2k-3)\b-4}{2\g-(k-3)\b-4}\prod_{i=1}^k
\frac{(2\g-(i-3)\b-2)(2\g-(i-3)\b-4)((\g-(i-3)\b-4)}{i\b((i-1)\b+2)(\g-(i-1)\b)}.$$

\smallskip
When $k$ is small enough, $Y'_k$ is an irreducible module whose highest weight 
is given by Proposition 6.1 below. But the formula above may give a nonzero 
integer when $k$ is too big for the hypothesis of this Proposition to hold. 
We check case by case that, nevertheless, this integer is still the dimension of an irreducible 
module, or possibly  {\it the opposite} of the dimension of an irreducible 
module. This means that $Y'_k$ should be interpreted as a virtual module, which
is a true module for small $k$, possibly the opposite
of  a module for intermediate values 
of $k$, and zero for     $k$
sufficiently large. In the second situation, we put a minus sign before 
the highest weight of the corresponding module in the lists below. 

\bigskip\noindent
{\bf a.}
$Y_k(\b)$ for $\fsp_{2l}$ (note that we have the fold of $\fsl_{2l+2}$):

$$\begin{array}{cccccccccccc}
{\scriptstyle k}&{\scriptstyle  0}&{\scriptstyle  1}&{\scriptstyle 2}&\ldots & {\scriptstyle l}&
{\scriptstyle l+1,l+2}&{\scriptstyle l+3}& \ldots & 
{\scriptstyle 2l+2}&{\scriptstyle 2l+3}&{\scriptstyle\ge 2l+4}\\
Y'_k & \BC&2\o_1&2\o_2&\cdots & 2\o_l& 0& -2\o_l&\cdots & -2\o_1&\BC&0
\end{array}$$

\bigskip\noindent
{\bf b.} $Y_k(\b)$ for $\fsl_{l+1}$, for $l=2m-1$ odd and $l=2m$ even, respectively: 

$$
\begin{array}{ccccccccccccccccc}
{\scriptstyle k}&{\scriptstyle  0}&{\scriptstyle  1}&\cdots &{\scriptstyle m-1}&
{\scriptstyle m}& & & & & \\
Y'_k & \CC & 
\o_1+\o_{2m-1} &\cdots & \o_{m-1}+\o_{m+1} & 2\o_m & & & & & \\ 
 & & & &  {\scriptstyle m+1}&
{\scriptstyle m+2}&\cdots & {\scriptstyle 2m+1}&
 {\scriptstyle 2m+2}& {\scriptstyle \ge 2m+3} \\
 & & & & -2\o_m & -(\o_{m-1}+\o_{m+1}) &\cdots & 
-(\o_1+\o_m) & -\BC &  0 \\ 
 & & & & & & & & \\
{\scriptstyle k}&{\scriptstyle  0}&{\scriptstyle  1}&\ldots &{\scriptstyle m}&
{\scriptstyle m+1} & & & \\
Y'_k & \CC & 
\o_1+\o_{2m} & \cdots & \o_{m}+\o_{m+1} & 0 & & & \\
 & & & & & {\scriptstyle m+2}&
 \ldots & {\scriptstyle 2m+1}&
 {\scriptstyle 2m+2}& {\scriptstyle \ge 2m+3} \\
 &  & & &  & -(\o_{m}+\o_{m+1})
& \cdots & -(\o_1+\o_{2m}) & -\BC &  0
\end{array}
$$

\bigskip\noindent
{\bf c.} $Y_k(\b)$ for $\fso_{2l+1}$, for $l=2m-1$ odd and $l=2m$ even, respectively:

$$\begin{array}{cccccccccc}
{\scriptstyle k}&{\scriptstyle 0}& {\scriptstyle 1}& \cdots &{\scriptstyle m-1}& 
{\scriptstyle m}&{\scriptstyle m+1}&\cdots & {\scriptstyle l=2m-1}& {\scriptstyle \ge l+1} \\
Y'_k & \CC & \o_2 &\cdots & \o_{2m-2} &2\o_{2m-1}& \o_{2m-3}&\cdots & \o_1 & 0\\
 & & & & & & & &  \\
{\scriptstyle k}&{\scriptstyle 0}& {\scriptstyle 1}& \cdots &{\scriptstyle m-1}& 
{\scriptstyle m}&{\scriptstyle m+1}&\cdots & {\scriptstyle l=2m}&  {\scriptstyle \ge l+1} \\
Y'_k & \CC & \o_2 &\cdots & \o_{2m-2} &2\o_{2m}& \o_{2m-1}&\cdots & \o_1 & 0
\end{array}
$$

\bigskip\noindent
{\bf d.}
$Y_k(\b)$ for $\fso_{2l}$, $l\geq 4 $, for $l=2m-1$ odd and $l=2m$ even, respectively:

$$\begin{array}{ccccccccccc}
{\scriptstyle k}&{\scriptstyle 0}& {\scriptstyle 1}& \cdots &{\scriptstyle m-2}& 
{\scriptstyle m-1,m} & {\scriptstyle m+1}&\cdots& {\scriptstyle 2m-1}& {\scriptstyle \ge 2m} \\
Y'_k & \CC & \o_2 &\cdots & \o_{2m-4}& \o_{2m-2}+\o_{2m-1}& \o_{2m-3}&\cdots & \o_1 & 0\\
 & & & & & & & & & \\
{\scriptstyle k}&{\scriptstyle 0}& {\scriptstyle 1}& \cdots &{\scriptstyle m-1}& 
{\scriptstyle m}&{\scriptstyle m+1}&\cdots & {\scriptstyle   2m}&  {\scriptstyle \ge l+1} \\
Y'_k & \CC & \o_2 &\cdots & \o_{2m-2} &0& \o_{2m-2}&\cdots &\CC & 0
\end{array}
$$

\bigskip\noindent
{\bf e.}
$Y_k(\b)$ for the exceptional Lie algebras:
 
$$\begin{array}{lccccccccc} 
 & {\scriptstyle k}&{\scriptstyle 0}& {\scriptstyle 1}& {\scriptstyle 2}& {\scriptstyle 3}&
{\scriptstyle 4}& {\scriptstyle 5}& {\scriptstyle 6}& {\scriptstyle \ge 7} \\   
\fg_2 \quad & Y'_k & \BC & \o_2 & 2\o_1 & \o_1 &  0 & 0 & 0 & 0 \\
 \ff_4 \quad & Y'_k& \BC & \o_1 & 2\o_4 & \o_3 & \o_4& 0 & 0  & 0 \\
\fe_6 \quad & Y'_k& \BC & \o_2 & \o_1+\o_6 & \o_1+\o_6& \o_2 & \BC & 0 & 0 \\
\fe_7 \quad & Y'_k& \BC & \o_1 & \o_6 & 2\o_7 &  0 & 0 & 0 & 0\\
\fe_8 \quad & Y'_k& \BC & \o_8 & \o_1 & 0& -\o_1& -\o_8& -\BC&  0
\end{array}
$$
 
\bigskip\noindent
{\bf f.} $Y_k(\g)$: 
$$\begin{array}{lccccc} 
 & {\scriptstyle k}&{\scriptstyle 2}& {\scriptstyle 3}& {\scriptstyle 4}& {\scriptstyle \ge 5}\\
\fsp_{2l}\quad & Y''_k & \o_2 & 0 & 0 & 0 \\
\fsl_l & Y''_k& \o_1+\o_{l-1} & \CC & 0 & 0 \\
\fso_{m}  & Y''_k& 2\o_1 & 0 & 0 & 0\\
{\mathfrak e},{\mathfrak f},{\mathfrak g} & Y''_k  & 0 & -\fg & -\CC & 0 
\end{array}
$$

\medskip
  Note that the $Y_k(\b)$'s travel nicely along the Dynkin diagram
in a wave (that gets reflected when it hits the end of a diagram
or collides with something else in the diagram, becoming negative if
there is no arrow).\medskip

\subsection{Gradings and the $Y_k(\b)$'s}
Let $\b_1=\a_0$ denote the highest root, let $\b_2=\ta$ denote a highest long root orthogonal
to $\a_0$, $\b_3$ a  highest long root orthogonal to $\a_0$ and $\b_2$, etc...

\begin{proposition}\label{prop61} If $\s_i=\b_1+\b_2+\b_3+\cdots +\b_i$ is dominant for $i\le k$, 
then $\s_k$ is the highest weight of an irreducible component of $S^k\fg$.
\end{proposition}

This module turns out to be $Y'_k$, the module that we identified from its dimension. 
Of course this does not explain what happens when $Y'_k$ is only a virtual module, and in fact 
there are also cases for which $Y'_k$ is an actual module whose presence is not 
accounted for by the Proposition. To be precise, this happens when $k\ge m$ 
for $\fso_{4m-3}$, $\fso_{4m-2}$, $\fso_{4m-1}$ or $\fso_{4m}$, when $k\ge 2$ 
for $\fg_2$, when $k\ge 3$ for $\ff_4$ and $\fe_6$. 
  
\smallskip
Note that the dominance condition in the hypothesis of   Proposition
\ref{prop61} is not automatic, 
and will be essential 
in the following construction. We associate to the $k$ roots $\b_1,\ldots ,\b_k$
a $\ZZ^k$-grading of $\fg$, 
$$\fg_{l_1\cdots l_k}=\{X\in\fg, [H_{\b_i},X]=l_iX, i=1...k\}.$$

\begin{lemma} Suppose that $l_1,\ldots ,l_k\ge 0$ and $\fg_{l_1\cdots l_k}\ne 0$. 
Then $l_1+\cdots +l_k\le 2$. 
\end{lemma}

\medskip
Note that for a given $\fg$ which is not $\fso_n$ for some $n\ge 5$, there is no ambiguity 
in defining the integer $\b$. This implies that for any $k$ as above, 
the components $\fg_{0..1..1..0}$ of our $k$-dimensional grading have the 
same dimension, $\b$ by Proposition \ref{prop51}, in particular they are nonzero.

If $\fg=\fso_n$ for some $n\ge 5$, the Lie subalgebra we denoted $\fh$ is the product of 
$\fsl_2$ and $\fso_{n-4}$, and we can choose for $\b_2$ the highest root of either algebra. 
If we choose that of $\fsl_2$, we cannot go further: $\b_1+\b_2+\b_3$ will not be 
dominant. If we choose the highest root of $\fso_{n-4}$, we can go further, but there is
no more choice, we can only take $\b_i=\e_{2i-1}+\e_{2i}$ and again the components 
$\fg_{0..1..1..0}$ of the grading have the same dimension, four. 
 
\exam For $\fg=\fe_7$, the highest root is $\a_0=\b_1=\o_1$. The highest root orthogonal 
to $\a_0$ is the highest root of a subsystem of type $D_6$. We get $\b_2=\o_6-\o_1$  
and $\b_1+\b_2=\o_6$ is dominant. For the next step, the roots orthogonal both to 
$\b_1$ and $\b_2$, i.e., both to $\o_1$ and $\o_6$, form a reducible subsystem of type
$D_4\times A_1$, and we have two candidates for the next highest root. If we choose,
$\b_3=\o_4-\o_1-\o_6$, the highest root of the $D_4$ part, then $\b_1+\b_2+\b_3$ is 
not dominant. The only possible choice is therefore $\b_3=\a_7=2\o_7-\o_6$, for 
which  $\b_1+\b_2+\b_3=2\o_7$ is dominant. Then the process stops. 

We obtain a three-dimensional grading of $\fe_7$ with three types of nonzero components:
the six components $\fg_{\pm 200}, \fg_{0\pm 20}, \fg_{00\pm 2}$ are one dimensional;
the twelve components $\fg_{\pm 1\pm 10}, \fg_{\pm 10\pm 1}, \fg_{0\pm 1\pm 1}$ have 
dimension eight and must be interpreted as copies of the (complexified) octonions;
the central component $\fg_{000}=\CC^3\op\fso_8$. This
is very close to the triality
construction of $\fe_7$ as $\fg(\OO,\HH)$ \cite{LMadv}.

\begin{corollary} 
For each $1\le i<k$, we have $(2\rho,\b_i-\b_{i+1})=\b(\a_0,\a_0)$. 
\end{corollary} 

\proof Recall that $2\rho$ is the sum of the positive roots. If $\g$ is such a 
positive root, then $\fg_{\g}$ is contained in one of the $\fg_{l_1...l_k}$, 
and the fact that $\b_1+\cdots +\b_j$ is dominant for all $j$ implies that 
$l_1+\cdots +l_j\ge 0$ for all $j$. Conversely, such a component $\fg_{l_1...l_k}$
is a sum of positive root spaces, except of course the central component $\fg_{0...0}$. 
 
The integer $\g(H_{\b_i})$ is equal to $2$ if $\fg_{\g}=\fg_{0...2...0}$ with the $2$ 
in position $i$, $1$ if $\fg_{\g}\subset\fg_{0..1..1..0}$ or $\fg_{0..1..-1..0}$, both
of dimension $\b$, $1$ again for $\fg_{0..1..0}$, of dimension with the $1$ in position $i$, 
and zero otherwise. Since $\fg_{0..1..0}$ has dimension $2\g-8-2(k-2)\b$, we conclude that 
$$2\rho(H_{\b_i})=(k-1+k-i)\b+2\g-8-2(k-2)\b+2=2\g+(3-i)\b-6,$$
and the claim follows.\qed  

\begin{corollary} 
The Casimir eigenvalue of $Y_k(\b)$ is $2kt-k(k-1)\b$.
\end{corollary} 

\noindent {\it Proof of the Lemma}.
Let $\th$ be some root such that $\fg_{\th}\subset \fg_{l_1\cdots l_k}$. 
If some $l_i$ equals $2$, then $\th$ must equal $\b_i$ and the other coefficients
vanish. So we suppose that $l_{i_1}=\cdots =l_{i_p}=1$, and the other coefficients
are zero. Using the orthogonality of the $\b_i$'s, we can write 
$$\th=\frac{1}{2}(\b_{i_1}+\b_{i_2}+\cdots +\b_{i_p}+\g),$$
where $\g$ is orthogonal to the $\b_i$'s. Suppose that $i_1$ is smaller
than the other $i_q$'s and apply the symmetry $s=s_{i_2}\cdots s_{i_p}$. 
We conclude that 
$$s(\th)=\frac{1}{2}(\b_{i_1}-\b_{i_2}-\cdots -\b_{i_p}+\g)$$
is again a root. Since $(s(\th),\b_1+\cdots +\b_{i_1})>0$, it must be a 
positive root. But
$$(\b_1+\cdots +\b_k,s(\th))=1-(p-1)=2-p.$$
Since $\b_1+\cdots +\b_k$ is supposed to be dominant, 
this must be a nonnegative integer. Hence $p\ge 2$, which is what we wanted to prove. \qed

\medskip\noindent  {\it Proof of Proposition 6.1.}
The case $k>3$ only happens for the classical Lie algebras, for 
which we can exhibit a highest weight vector of weight $\s_k$ as a determinant, or a 
Pfaffian, in terms of the basis of the natural representation preserved by the 
maximal torus:
\begin{eqnarray}
\nonumber
\fsl_n\qquad &   \det (e_p\ot e_{n+1-q}^*)_{1\le p,q\le k}, \\
\nonumber
\fsp_{2n}\qquad &   \det (e_p\circ e_q)_{1\le p,q\le k}, \\
\nonumber
\fso_n\qquad &   {\rm Pf} (e_p\wedge e_q)_{1\le p,q\le k}.
\end{eqnarray}

\smallskip So we focus on the case $k=3$. We have a $\ZZ^3$-grading of $\fg$ 
with $6$ terms of dimension one, $12$ of dimension $\b$, $6$ of dimension 
$2\g-2\b-8$, and the central term $\fg_{000}$. The six terms of type $\fg_{200}$
can be represented as the vertices of a square pyramid; then the twelve terms 
of type $\fg_{110}$ are the middle points of the edges. 

\begin{center}
\setlength{\unitlength}{3mm}
\begin{picture}(20,25)(-9,-12)
\put(-.2,-.2){$\circ$}
\put(1.5,11.5){$\fg_{002}$}
\put(12,1){$\fg_{200}$}
\put(6.8,6.2){$\fg_{101}$}

\put(.2,5.4){$\bullet$}
\put(.2,-1.55){$\bullet$}
\put(5.85,-.4){$\bullet$}
\put(-5.5,-.15){$\bullet$}
\put(-.75,1.15){$\bullet$}
\put(-.7,-5.85){$\bullet$}

\put(.7,-1.3){\line(4,1){5.3}}
\put(.3,-1.25){\line(-4,1){5.35}}
\put(-5.3,.05){\line(1,1){5.6}}
\put(6.1,.1){\line(-1,1){5.5}}
\put(.55,-1.15){\line(0,1){6.7}}

\put(-5.3,.15){\line(4,1){4.7}}
\put(-.4,1.45){\line(1,5){.85}}
\put(-.2,1.3){\line(5,-1){6.2}}
\put(-.35,-5.55){\line(1,5){.85}}
\put(-.4,-5.35){\line(0,1){6.8}}
\put(-.2,-5.35){\line(6,5){6.2}}
\put(-.55,-5.4){\line(-5,6){4.55}}

\put(.5,10.9){$\bullet$}
\put(.55,-3.15){$\bullet$}
\put(11.75,-.45){$\bullet$}
\put(-10.8,-.2){$\bullet$}
\put(-1.3,2.3){$\bullet$}
\put(-1.3,-11.5){$\bullet$}

\thicklines
\put(1.1,-2.8){\line(4,1){10.7}}
\put(.6,-2.7){\line(-4,1){10.9}}
\put(-10.5,.1){\line(1,1){11.1}}
\put(11.9,.1){\line(-1,1){11}}
\put(.9,-2.5){\line(0,1){13.5}}

\put(-10.5,.2){\line(4,1){9.3}}
\put(-.85,2.8){\line(1,5){1.7}}
\put(-.7,2.55){\line(5,-1){12.5}}
\put(-.9,-11){\line(1,5){1.6}}
\put(-1,-11){\line(0,1){13.6}}
\put(-.8,-11){\line(6,5){12.8}}
\put(-1.25,-11.05){\line(-5,6){9.15}}

\put(.6,3.9){$\ast$}
\put(-.35,6.6){$\ast$}
\put(-5.1,5.4){$\ast$}
\put(6.2,5.2){$\ast$}
\put(5.25,.9){$\ast$}
\put(6.2,-1.8){$\ast$}
\put(-5.1,-1.6){$\ast$}
\put(-6.05,1.1){$\ast$}
\put(-.35,-7.3){$\ast$}
\put(-1.3,-4.6){$\ast$}
\put(-6.05,-5.8){$\ast$}
\put(5.25,-6){$\ast$}

\end{picture}
\end{center}
\medskip
We have chosen vectors $X_{\pm\b_i}$ generating three commuting $\fsl_2$-triples. 
We next choose generators $X_{\a}$ for the roots $\a\in\Phi_{110}$ such that the 
corresponding root  space $\fg_{\a}\subset\fg_{110}$. We then  
take bases  of $\fg_{-110}$, $\fg_{1-10}$ and $\fg_{-1-10}$ by letting 
$$X_{\a-\b_1}=[X_{-\b_1},X_{\a}], \quad X_{\a-\b_2}=[X_{-\b_1},X_{\a}], \quad 
X_{\a-\b_1-\b_2}=[X_{-\b_1},X_{\a-\b_2}].$$
An easy consequence of the Jacobi identity is that $[X_{\b_1},X_{\a-\b_1}]
=X_{\a}$.  
Note that $\a-\b_1-\b_2=s_{\b_1}s_{\b_2}(\a)$ is a root, and that its opposite
$$
\a'=\b_1+\b_2-\a
$$ corresponds to another root space in $\fg_{110}$.
(We will use the notation $\th'=\b_1+\b_2-\th$ repeatedly in what follows.)
 Then $[X_{\a},X_{\a'-\b_2}]$, 
which is equal to $[X_{\a'},X_{\a-\b_2}]$ by the Jacobi identity, is a nonzero 
multiple of $X_{\b_1}$. We normalize our root vectors from $\fg_{110}$ so that this
multiple is in fact $X_{\b_1}$ itself. This means that  for all $\a\in\Phi_{110}$,
 $$K(X_{\b_1},X_{-\b_1})=
K([X_{\a},X_{\a'-\b_2}],X_{-\b_1})=K(X_{\a},[X_{\a'-\b_2},X_{-\b_1}])=-
K(X_{\a},X_{-\a}).$$

We use the same normalization for $\fg_{101}$ and $\fg_{011}$. 
Note that $2K(X_{\b_1},X_{-\b_1})=K(H_{\b_1},H_{\b_1})=2/(\b_1,\b_1)$, twice the inverse 
of the square length of a long root. In particular, $K(X_{\a},X_{-\a})$ does not depend
on $\a\in\Phi_{\pm 1\pm 1 0}$, $\Phi_{\pm 10\pm 1}$ or $\Phi_{0\pm 1\pm 1}$.
 
\medskip Now we introduce the symmetric tensor 
$$S_{12}=\sum_{\a\in\Phi_{110}}X_{\a}X_{\a'}.$$
Here and in what follows, if $X,Y\in \fg$ or a symmetric power of $\fg$,
$XY$ will denote the symmetric product $X\circ Y$.

\begin{lemma}\label{lem65} For $Y,Z\in\fg_{-1-10}$, the bilinear forms 
$$\sum_{\a\in\Phi_{110}}K(X_{\a},Y)K(X_{\a'},Z)\quad and\quad
K([X_{\b_1},Y],[X_{\b_2},Z])$$
are multiples of one another. 

In particular, $S_{12}$ is invariant 
under the common centralizer of $\b_1$ and $\b_2$.
\end{lemma}
 
\proof Let $Y=X_{-\b'}$ and $Z=X_{-\g'}$ for some roots $\b,\g\in\Phi_{110}$. 
Then $[X_{\b_1},Y]=X_{\b-\b_2}$ and  $[X_{\b_2},Z]=X_{\g-\b_1}$, hence
$$K([X_{\b_1},Y],[X_{\b_2},Z])=K(X_{\b-\b_2},X_{\g-\b_1})=-K(X_{-\b'},X_{\g})=
\delta_{\b',\g}K(X_{-\b_1},X_{\b_1}).$$
Since $\sum_{\a\in\Phi_{110}}K(X_{\a},Y)K(X_{\a'},Z)=\delta_{\b',\g}K(X_{-\b_1},X_{\b_1})^2$,
the claim follows. \qed

\medskip We deduce a different proof of Proposition 3.1. We must prove that $S^2\fg$ 
contains a tensor of weight $\a_0+\ta=\b_1+\b_2$ which is a highest weight vector, i.e., which
is annihilated by any positive root vector. 

\begin{corollary}
The tensor $\Sigma=X_{\b_1}X_{\b_2}-\frac{1}{2}S_{12}\in S^2\fg$ 
is a highest weight vector of weight $\b_1+\b_2$. 
\end{corollary}

\proof We must prove that $\Sigma$ is annihilated by any positive root vector. 
Since $\b_1$ and $\b_1+\b_2$ are both dominant, a positive root must belong either
to $\Phi_{00}$, $\Phi_{1-1}$ or $\Phi_{pq}$ with $p,q\ge 0$ and $p+q>0$. Since
$\fg_{1+p,1+q}=0$, our assertion is clear for the latter case.  For the first case, 
it follows from the previous lemma. 

There remains to prove that $ad(X_{\theta-\b_2})\Sigma =0$ for $\theta\in\Phi_{11}$.
We use the structure constants $N_{\mu,\nu}$ such that 
$[X_{\mu}, X_{\nu}]=N_{\mu,\nu}X_{\mu+\nu}$. 
Note that $ad(X_{\theta-\b_2})X_{\b_1}X_{\b_2} =X_{\b_1}X_{\theta}$.
 
 On the other hand,
$$ad(X_{\theta-\b_2})S_{12} =
2\sum_{\a\in\Phi_{11}}N_{\theta-\b_2,\a}X_{\a+\theta-\b_2}X_{\a'}.$$
But $\a+\theta-\b_2$ belongs to $\Phi_{200}$, hence must be equal to $\b_1$, 
which forces $\a=\theta'$. Our normalization is $N_{\theta-\b-2,\theta'}=-1$, 
thus $ad(X_{\theta-\b_2})S_{12} =2X_{\b_1}X_{\theta}$. This concludes the proof.\qed

\medskip

Now we define a tensor $T\in \fg_{110}\ot \fg_{101}\ot \fg_{011}\subset S^3\fg$, with the 
help of which we will construct a highest weight vector of weight $\b_1+\b_2+\b_3$:
\begin{eqnarray}
 T=\sum_{\substack{\a\in\Phi_{011},\b\in\Phi_{101},\g\in\Phi_{110}\\ 
\a+\b+\g=\b_1+\b_2+\b_3}}
N_{\b_3-\a,\b_1-\b}X_{\a}X_{\b}X_{\g}.
\end{eqnarray}

We will need the following properties of the structure constants. 

\begin{lemma}
If $\mu\in\Phi_{1-10}$ and $\nu\in\Phi_{01-1}$, then
\begin{eqnarray}
N_{\mu,\nu} &= &-N_{\mu+\b_2,\nu-\b_2}, \\
N_{\mu,\nu} &= &N_{\mu-\b_1,\nu}, \\
N_{\mu,\nu} &= &N_{\nu,-\mu-\nu}.
\end{eqnarray}
\end{lemma}

Of course, we have similar identities when we permute the indices, 
e.g., $N_{\mu,\nu} = -N_{\mu+\b_3,\nu-\b_3}$ if $\mu\in\Phi_{10-1}$ 
and $\nu\in\Phi_{0-11}$.

\proof By definition, $X_{\mu}=[X_{-\b_2},X_{\mu+\b_2}]$. In the Jacobi identity
$$[[X_{-\b_2},X_{\mu+\b_2}],X_{\nu}]+[[X_{\mu+\b_2},X_{\nu}],X_{-\b_2}]+
[[X_{\nu},X_{-\b_2}],X_{\mu+\b_2}]=0,$$
the first bracket of the second term is in $\fg_{12-1}$, hence equal to zero. 
Since $[X_{\nu},X_{-\b_2}]=-X_{\nu-\b_2}$, the first identity follows. 
To prove the second one, we use the Jacobi identity
$$[[X_{-\b_1},X_{\mu}],X_{\nu}]+[[X_{\mu},X_{\nu}],X_{-\b_1}]+
[[X_{\nu},X_{-\b_1}],X_{-\mu}]=0.$$
Here the first bracket of the last term is in $\fg_{-11-1}$, hence equal to zero
and we deduce the second identity. 
Finally, the invariance of the Killing form gives
$$N_{\mu,\nu}K(X_{\mu+\nu},X_{-\mu-\nu})=K([X_{\mu},X_{\nu}],X_{-\mu-\nu})=
K(X_{\mu},[X_{\nu},X_{-\mu-\nu}])=N_{\nu,-\mu-\nu}K(X_{\mu},X_{\mu}).$$
But in the normalization we use, we have seen that $K(X_{\mu+\nu},X_{-\mu-\nu})=
K(X_{\mu},X_{\mu})$, and the third identity follows. \qed

\begin{proposition}
The tensor $\Theta\in S^3\fg$ defined as
$$\Theta=X_{\b_1}X_{\b_2}X_{\b_3}-X_{\b_1}S_{23}-X_{\b_2}S_{13}-X_{\b_3}S_{12}+T,$$
is a highest weight vector of weight $\b_1+\b_2+\b_3$.
\end{proposition}

\proof We must check that $\Theta$ is annihilated by any positive root vector $X_{\mu}$. 
If $\mu(H_{\b_i})\ge 0$ for $i=1,2,3$, and at least one is positive, this is clear 
since $X_{\mu}$ annihilates every space of type $\fg_{200}$  or $\fg_{110}$. 
If these three integers vanish, that is, $X_{\mu}\in\fg_{000}$, this follows 
from the fact that for $X\in\fg_{0-1-1}$, $Y\in\fg_{-10-1}$ and $Z\in\fg_{-1-10}$,
$$\sum_{\substack{\a\in\Phi_{011},\b\in\Phi_{101},\g\in\Phi_{110}\\ 
\a+\b+\g=\b_1+\b_2+\b_3}}
N_{\a-\b_3,\b-\b_1}K(X_{\a},X)K(X_{\b},Y)K(X_{\g},Z)=
K([[X_{\b_1},Y],[X_{\b_2},Z]],[X_{\b_3},X]),$$
which shows that $T$ must be annihilated by any vector commuting with 
$X_{\b_1}$, $X_{\b_2}$ and  $X_{\b_3}$ -  and $X_{\mu}$ has this property.

Now, since $\mu$ is positive, we know that $\mu(H_{\b_1})$, $\mu(H_{\b_1})+\mu(H_{\b_2})$
and $\mu(H_{\b_1})+\mu(H_{\b_2})+\mu(H_{\b_3})$ are nonnegative, so if one of the 
$\mu(H_{\b_i})$'s is negative, $X_{\mu}$ must belong to $\fg_{1-10}$, $\fg_{10-1}$
or $\fg_{01-1}$. Since $[ \fg_{1-10},\fg_{01-1}]=\fg_{10-1}$, what remains to check
is that $\Theta$ is annihilated by any $Z\in\fg_{1-10}$ or $Y\in\fg_{01-1}$.
This is equivalent to the four identities
\begin{eqnarray}
\; [Z_3,T] &= (ad(Z)S_{23}) \circ X_{\b_1}, \\
\; [Y_1,T] &= (ad(Y) S_{13}) \circ  X_{\b_2}, \\ 
\; [Z_1,T] &= S_{13}\circ ad(Z)X_{\b_2} , \\
\; [Y_2,T] &= S_{12}\circ ad(Y) X_{\b_3} , 
\end{eqnarray}
where $[Z_3,T]$, for example, means that we take the bracket of $Z$ only with the
terms in $T$ coming from $\fg_{110}$.  

\medskip\noindent {\it Proof of (5)}. 
To prove the first identity, we can let $Z=X_{\d-\b_2}$ for some $\d\in\Phi_{110}$. 
Then
$$[Z_3,T]=\sum_{\substack{\a\in\Phi_{011},\b\in\Phi_{101}, \\ 
\a+\b=\d+\b_3}}N_{\b_3-\a,\b_1-\b}X_{\b_1}X_{\a}X_{\a'+\d-\b_2}.$$
But 
$$\begin{array}{rcll}
N_{\b_3-\a,\b_1-\b} & = & -N_{-\a,\b_3+\b_1-\b} & \qquad (2) \\
  & = & N_{\b_2-\a,\b_3-\b} & \qquad (3)\;  {\rm twice} \\
  & = & N_{\b_2-\a,\a-\d} &  \\
  & = & N_{\b_2-\a,\d-\b_2} &  \qquad (4) \\
  & = & N_{\a',\d-\b_2} & \qquad (3) \\
\end{array}$$
This allows us to write $[Z_3,T]$ as 
$$\sum_{\a\in\Phi_{011}}N_{\a',\d-\b_2}X_{\b_1}X_{\a}X_{\a'+\d-\b_2}
=\frac{1}{2}(ad(Z)(\sum_{\a\in\Phi_{011}}X_{\a}X_{\a'}))X_{\b_1}=(ad(Z)S_{23})\circ X_{\b_1}.$$
This proves (5). The proof of (6) is similar and will be left to the reader. 

\medskip\noindent {\it Proof of (8)}.
The proofs of (7) and (8) involve the same type of arguments and we will
focus on (8). 
We use the invariance of $S_{12}$  from Lemma \ref{lem65}. Let $Y=X_{\theta-\b_3}$, with 
$\theta\in\Phi_{011}$. We have $ad(X_{\theta})S_{12}=0$ since $\fg_{121}=0$. 
Thus for $\a\in\Phi_{011}$, 
$$ad(X_{\theta})ad(X_{-\a})S_{12}=ad([X_{\theta},X_{-\a}])S_{12}.$$
If $\a\ne\theta$, $[X_{-\theta},X_{\a}]$ is either zero, or a root vector in 
$\fg_{000}$, thus annihilating $S_{12}$. Hence
$$\begin{array}{rcl}
X_{\theta} ad([X_{\theta},X_{-\theta}])S_{12} & = &
\sum_{\a\in\Phi_{011}}X_{\a}ad(X_{\theta})ad(X_{-\a})S_{12} \\ & = & 
\sum_{\a\in\Phi_{011},\g\in\Phi_{110}}N_{-\a,\g'}N_{\g'-\a,\theta}X_{\a}X_{\g'-\a+\theta}X_{\g}\\
 & = & \sum_{\substack{\a\in\Phi_{011},\b\in\Phi_{101},\g\in\Phi_{110}\\ 
\a+\b+\g=\b_1+\b_2+\b_3}}N_{-\a,\g'}N_{\g'-\a,\theta}X_{\a}X_{\b+\theta-\b_3}X_{\g}.
\end{array}$$
But $N_{\g'-\a,\theta}=N_{\b-\b_3,\theta}=-N_{\b,\theta-\b_3}$ by (2), and 
$N_{-\a,\g'}=N_{-\a,\a+\b-\b_3}=N_{\b_3-\a,\a+\b-\b_3}=-N_{\b_3-\a,-\b}=-N_{\b_3-\a,\b_1-\b}$, 
where we used successively (3), (4) and (3) again. Thus 
$$\begin{array}{rcl}
X_{\theta} ad([X_{\theta},X_{-\theta}])S_{12} 
 & = & \sum_{\substack{\a\in\Phi_{011},\b\in\Phi_{101},\g\in\Phi_{110}\\ 
\a+\b+\g=\b_1+\b_2+\b_3}}N_{\b_3-\a,\b_1-\b}N_{\b,\theta-\b_3}X_{\a}X_{\b+\theta-\b_3}X_{\g}, \\
 & = & \sum_{\substack{\a\in\Phi_{011},\b\in\Phi_{101},\g\in\Phi_{110}\\ 
\a+\b+\g=\b_1+\b_2+\b_3}}N_{\b_3-\a,\b_1-\b}X_{\a}[X_{\b},X_{\theta-\b_3}]X_{\g} \\
 & = & -[Y_2,T].
\end{array}$$
There remains to compute $ad([X_{\theta},X_{-\theta}])S_{12}$. We have 
$[X_{\theta}, X_{-\theta}]=t_{\theta}H_{\theta}$, where $t_{\theta}$
can be computed as follows:  
$$t_{\theta}K(H_{\theta},H_{\theta})=K([X_{\theta}, X_{-\theta}],H_{\theta})
=K(X_{\theta}, [X_{-\theta},H_{\theta}])=2K(X_{\theta}, X_{-\theta}).$$
And since we know that $2K(X_{\theta}, X_{-\theta})=-2K(X_{\b_2}, X_{-\b_2})
=-K(H_{\b_2},H_{\b_2})$, we get that 
$$t_{\theta}=-\frac{K(H_{\b_2},H_{\b_2})}{K(H_{\theta},H_{\theta})}=
-\frac{(\theta,\theta)}{(\b_2,\b_2)}=-\frac{(\theta,\theta)}{(\b_1,\b_1)}.$$
Then $ad([X_{\theta},X_{-\theta}])S_{12}=t_{\theta}ad(H_{\theta})S_{12}$ is equal to
$$t_{\theta}\sum_{\g\in\Phi_{110}}(\g(H_{\theta})+\g'(H_{\theta}))X_{\g}X_{\g'}
=t_{\theta}(\b_1(H_{\theta})+\b_2(H_{\theta}))S_{12}.$$
But since $\theta\in\Phi_{011}$, $\b_1(H_{\theta})=\theta(H_{\b_1})=0$, while 
$$\b_2(H_{\theta})=\frac{(\b_2,\b_2)}{(\theta,\theta)}\theta(H_{\b_2})=-t_{\theta}^{-1}.$$
 We thus  get that $[Y_2,T]=S_{12}[Y,X_{\b_2}]$, as required. This concludes the proof 
of the identity (8), hence of Proposition 6.1.   \qed
\medskip

\subsection{Geometric interpretation of the $Y_k(\b)$'s}

Zak defines the Scorza varieties to be the 
smooth nondegenerate varieties extremal
for higher secant defects in the sense the the
defect of the $i$-th secant variety of $X^n\subset\BP^N$ is $i$ times the defect $\delta$ of
the first and that the $[\frac n\delta]$-th secant variety fills
the ambient space. He then goes on to classify the Scorza varieties,
all of which turn out to be homogeneous \cite{zak}.

More precisely, the Scorza varieties are given by the
projectivization of the rank one elements in
the  Jordan algebras
$\cJ_r(\BA)$, where  $\BA$ is the complexification of the reals, complex
numbers or quaternions, or, when $r=3$, the octonions \cite{chaput2}.

Recall that for a simple Lie algebra, the adjoint group has a unique closed 
orbit in $\PP\fg$, the projectivization $X_{ad}$ of the minimal nilpotent orbit.
This {\it  adjoint variety} parametrizes the highest root spaces in $\fg$. 

\begin{proposition}\label{prop69}
Let $Z_k \subset X_{ad}$ denote the shadow of a point 
of the closed orbit 
  $X_k(\b)\subset\BP Y_k(\b)$. Then the $Z_k $'s are   Scorza varieties
on the adjoint variety.
\end{proposition}

The shadow of a point is defined as follows: let $G$ denote the adjoint
algebraic group associated to the Lie algebra $\fg$. We have maps
$$X_k(\b)=G/Q \stackrel{q}{\longleftarrow} G \stackrel{p}{\longrightarrow} G/P=X_{ad},$$
and the shadow of $x\in X_k(\b)$ is the subvariety $p(q^{-1}(x))$ of $X_{ad}$. 
Tits \cite{tits} showed how determine these shadows   using   Dynkin diagrams,   
  and   Proposition \ref{prop69} follows from a straightforward case by case check. 

\exam Let $\fg$ be $\fsl_{l+1}$ or $\fso_{2l}$. On the Dynkin diagram of $\fg$, 
we let the $\ast$'s encode the highest weight of the fundamental representation, 
and the $\bullet$'s encode that of $Y'_k$. When we suppress the $\bullet$'s, we get
weighted diagrams encoding homogeneous varieties, respectively $\PP^{k-1}\times\PP^{k-1}$
and a Grassmannian $G(2,2k)$ which are two examples of Scorza varieties. 

\begin{center}
\setlength{\unitlength}{3mm}
\begin{picture}(20,7)(9,1.5)
\multiput(4,5)(2,0){7}{$\circ$}
\multiput(4.5,5.4)(2,0){6}{\line(1,0){1.6}} 
\put(8,5){$\bullet$}
\put(12,5){$\bullet$}
\put(4,3.5){\line(0,1){.5}}
\put(6.5,3.5){\line(0,1){.5}}
\put(4,3.5){\line(1,0){2.5}}
\put(4.8,2){$\PP^2$}
\put(14,3.5){\line(0,1){.5}}
\put(16.5,3.5){\line(0,1){.5}}
\put(14,3.5){\line(1,0){2.5}}
\put(14.8,2){$\PP^2$}
\multiput(4,5)(12,0){2}{$\ast$}

\multiput(23,5)(2,0){6}{$\circ$}
\multiput(23.5,5.4)(2,0){5}{\line(1,0){1.6}} 
\put(33.4,5.4){\line(1,1){1.2}} 
\put(33.4,5.2){\line(1,-1){1.2}}
\put(34.6,6.5){$\circ$}  
\put(34.6,3.5){$\circ$}  
\put(33,5){$\bullet$}
\put(23,3.5){\line(0,1){.5}}
\put(31.5,3.5){\line(0,1){.5}}
\put(23,3.5){\line(1,0){8.5}}
\put(25.5,2){$G(2,6)$}
\put(25,5){$\ast$}

\end{picture}
\end{center}

It is natural that the Scorza varieties arrive as subsets of polynomials
of degree $k$ on $\fg$ because the determinant on $\cJ_k(\BA)$ is
a polynomial of degree $k$. If we take the linear span of $Z_k$ and
then take the cone over the
degree $k$ hypersurface in $\langle Z_k\rangle$ with vertex a Killing-complement to
$\langle Z_k\rangle$, we obtain a hypersurface of degree $k$ in $\fg$.
$X_k(\b)$ parametrizes this space of hypersurfaces and its span gives
the space $Y_k(\b)$.\medskip

\subsection{Universal dimension formulas}
We finally extend our formula for the dimension of the Cartan powers
of $\fg$ to obtain a universal formula for the Cartan powers of the $Y_l(\b)$'s. 
Again, our approach is based on Weyl's dimension formula: we check that the 
relevant integers can be organized into strings whose extremities depend only
on Vogel's parameters $\b$ and $\g$. In fact, this really makes sense only in
type $A,B,D$, and in the exceptional cases (excluding $\ff_4$)  when $l=2$. 
In type $C$ and $F_4$, there are some strange compensations involving half 
integers, but the final formula holds in all cases. 

We will just give the main statements leading to   Theorem \ref{bigformula} below. 
Let $\s_l=\b_1+\cdots +\b_l$ denote the highest weight of $Y_l(\b)$. 
Let $\Phi_{l,i}$ denote the set of positive roots $\g$ such that $\s_l(H_{\g})=i$. 
By Lemma 6.2, we have 
$$\begin{array}{rcl}
\#\Phi_{l,1} & = & 2l(\g+2\b-4-\b l), \\
\#\Phi_{l,2} & = & \b\frac{l(l-1)}{2}+l, \\
\#\Phi_{l,i} & = & 0 \qquad {\rm for}\; i>2. 
\end{array}$$

\noindent {\it Facts.}
\begin{enumerate}
\item The values of $\rho(H_{\g})$, for $\g$ in $\Phi_{l,2}$, can be organized into 
$l$ intervals $[\g-(2l-i-3)\b/2-3,\g-(i-3)\b/2-3]$, where $1\le i\le l$. 
\item The values of $\rho(H_{\g})$, for $\g$ in $\Phi_{l,1}$, can be organized into 
$3l$ intervals $[\g/2-(i-1)\b/2,\g/2-(l-i-2)\b/2-3]$, $[i\b/2,\g-(l+i-4)\b/2-3]$
and $[(i-1)\b/2+1,\g-(l+i-3)\b/2-4]$, with  $1\le i\le l$. 
\end{enumerate}

\medskip
Applying the Weyl  dimension formula, we obtain:

\begin{theorem} \label{bigformula}
$$\dim (Y_l(\b))^{(k)}
=\prod_{i=1}^l\frac{\binom{2k+\g-(i-3)\b/2-3}{2k}\binom{k+\g-(l+i-3)\b/2-3}{k}
\binom{k+\g-(l+i-4)\b/2-4}{k}\binom{k+\g/2-(l-i-2)\b/2-3}{k}}
{\binom{2k+\g-(2l-i-3)\b/2-4}{2k}\binom{k+i\b/2-1}{k}
\binom{k+(i-1)\b/2}{k}\binom{k+\g/2-(i-1)\b/2-1}{k}}.$$
\end{theorem}

\vspace{1cm}

{\small
\noindent {\sc Joseph M. Landsberg}, 
Department of Mathematics,
  Texas A\&M University,
  Mailstop 3368,
  College Station, TX 77843-3368, USA

\noindent {\rm E-mail}: jml@math.tamu.edu 

\medskip\noindent
{\sc Laurent Manivel},  Institut Fourier, UMR 5582 du CNRS 
Universit\'e Grenoble I, BP 74, 38402 Saint Martin d'H\`eres cedex, FRANCE

\noindent {\rm E-mail}: Laurent.Manivel@ujf-grenoble.fr

\end{document}